\documentclass[a4paper,12pt,twoside]{article}
\usepackage[dvips]{graphics}
\usepackage{epsf}
\usepackage{graphicx}
\usepackage{epsfig}
\usepackage{amsmath}
\usepackage[psamsfonts]{amssymb}
\usepackage[psamsfonts]{eucal}
\usepackage{amsthm}
\usepackage[T1]{fontenc}
\usepackage{textcomp}
\usepackage{mathptmx}
\usepackage[scaled]{helvet}
\usepackage{perso}

\swapnumbers
\newtheorem{defi}[subsection]{Definition}
\newtheorem{defis}[subsection]{Definitions}
\newtheorem{theo}[subsection]{Theorem}

\newtheorem{coro}[subsection]{Corollary}

\newcommand{\field}[1]{\mathbb{#1}}

\newcommand{\RR}{\field{R}}

\def\id{\mathop{\rm id}\nolimits}

\def\orth{\mathop{\rm orth}\nolimits}
\def\toto{\mathop{\rightrightarrows}\limits}

\title{Lie, symplectic and Poisson groupoids\\
        and their Lie algebroids}
\author{Charles-Michel Marle\\
  Universit{\'e} Pierre et Marie Curie\\
  Paris,  France}

\begin{document}
\maketitle

\section*{Introduction}

Groupoids are mathematical structures able to describe symmetry
properties more general than those described by groups. They were
introduced (and named) by H.~Brandt in 1926. Around 1950,
Charles~Ehresmann used groupoids with additional structures
(topological and differentiable) as essential tools in topology
and differential geometry. In recent years, Mickael Karasev,  Alan
Weinstein and Stanis\l aw Zakrzewski independently discovered that
symplectic groupoids can be used for the construction of
noncommutative deformations of the algebra of smooth functions on
a manifold, with potential applications to quantization. Poisson
groupoids were introduced by Alan Weinstein as generalizations of
both Poisson Lie groups and symplectic groupoids.

We present here the main definitions and first properties relative
to groupoids, Lie groupoids, Lie  algebroids, symplectic and
Poisson groupoids and their Lie algebroids.

\section{Groupoids}

\subsection{What is a groupoid?}

Before stating the formal definition  of a groupoid, let us
explain, in an informal way, why it is a very natural concept. The
easiest way to understand that concept is to think of two sets,
$\Gamma$ and $\Gamma_0$. The first one, $\Gamma$, is called the
\emph{set of arrows} or \emph{total space} of the groupoid, and
the other one, $\Gamma_0$, the \emph{set of objects} or \emph{set
of units} of the groupoid. One may think of an element $x\in
\Gamma$ as an arrow going from an object (a point in $\Gamma_0$)
to another object (another point in $\Gamma_0$). The word \lq\lq
arrow\rq\rq\ is used here in a very general sense: it means a way
for going from a point in $\Gamma_0$ to another point in
$\Gamma_0$. One should not think of an arrow as a line drawn in
the set $\Gamma_0$ joining the starting point of the arrow to its
end point: this happens only for some special groupoids. Rather,
one should think of an arrow as living outside $\Gamma_0$, with
only its starting point and its end point in $\Gamma_0$, as shown
on Figure~1.
\par\smallskip
The following ingredients enter the definition of a groupoid.

\begin{figure}[h]
\includegraphics{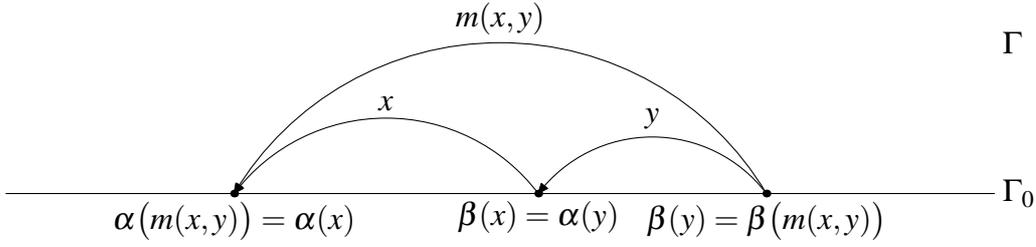}
\caption{Two arrows $x$ and $y\in\Gamma$, with the target of $y$,
$\alpha(y)\in \Gamma_0$, equal to the source of $x$,
$\beta(x)\in\Gamma_0$, and the composed arrow $m(x,y)$.}
\end{figure}

\begin{description}

\item{--} Two maps $\alpha:\Gamma\to\Gamma_0$ and $\beta:\Gamma\to
\Gamma_0$, called the {\it target map\/} and the {\it source
map\/} of the groupoid. If $x\in\Gamma$ is an arrow,
$\alpha(x)\in\Gamma_0$ is its end point and $\beta(x)\in\Gamma_0$
its starting point.

\item{--} A {\it composition law\/} on the set of arrows; we can
compose an arrow $y$ with another arrow $x$, and get an arrow
$m(x,y)$, by following first the arrow $y$, then the arrow $x$. Of
course, $m(x,y)$ is defined if and only if the target of $y$ is
equal to the source of $x$. The source of $m(x,y)$ is equal to the
source of $y$, and its target is equal to the target of $x$, as
illustrated on Figure~1. It is only by convention that we write
$m(x,y)$ rather than $m(y,x)$: the arrow which is followed first
is on the right, by analogy with the usual notation $f\circ g$ for
the composition of two maps $g$ and $f$. When there will be no
risk of confusion we will write $x\circ y$, or $x.y$, or even
simply $xy$ for $m(x,y)$.The composition of arrows is associative.

\item{--} An {\it embedding\/} $\varepsilon$ of the set $\Gamma_0$
into the set $\Gamma$, which associates a unit arrow
$\varepsilon(u)$ with each $u\in\Gamma_0$. That unit arrow is such
that both its source and its target are $u$, and it plays the role
of a unit when composed with another arrow, either on the right or
on the left: for any arrow $x$,
$m\Bigl(\varepsilon\bigl(\alpha(x)\bigr),x\Bigr)=x$, and
$m\Bigl(x, \varepsilon\bigl(\beta(x)\bigr)\Bigr)=x$.

\item{--} Finally, an {\it inverse map\/} $\iota$ from the set of
arrows onto itself. If $x\in\Gamma$ is an arrow, one may think of
$\iota(x)$ as the arrow $x$ followed in the reverse sense. We will
often write $x^{-1}$ for $\iota(x)$.
\end{description}

\par\smallskip
Now we are ready to state the formal definition of a groupoid.

\begin{defi}
A \emph{groupoid} is a pair of sets $(\Gamma,\Gamma_0)$ equipped
with the structure defined by the following data:
\begin{description}

\item{--} an injective map $\varepsilon:\Gamma_0\to\Gamma$, called
the \emph{unit section} of the groupoid;

\item{--} two maps $\alpha:\Gamma\to\Gamma_0$ and
$\beta:\Gamma\to\Gamma_0$, called, respectively, the \emph{target
map} and the \emph{source map\/}; they satisfy
\begin{equation}
\alpha\circ\varepsilon=\beta\circ\varepsilon=\id_{\Gamma_0}\,;
\end{equation}

\item{--} a composition law $m:\Gamma_2\to \Gamma$, called the
\emph{product}, defined on the subset $\Gamma_2$ of
$\Gamma\times\Gamma$, called the \emph{set of composable
elements},
\begin{equation}
\Gamma_2=\bigl\{\,(x,y)\in\Gamma\times\Gamma;
\beta(x)=\alpha(y)\,\bigr\}\,,
\end{equation}
which is associative, in the sense that whenever one side of the
equality
 \begin{equation}
 m\bigl(x,m(y,z)\bigr)=m\bigl(m(x,y),z\bigr)
 \end{equation}
is defined, the other side is defined too, and the equality holds;
moreover, the composition law $m$ is such that for each
$x\in\Gamma$,
\begin{equation}
 m\Bigl(\varepsilon\bigl(\alpha(x)\bigr),x\Bigr)
 =m\Bigl(x,\varepsilon\bigl(\beta(x)\bigr)\Bigr)=x\,;
\end{equation}

\item{--} a map $\iota:\Gamma\to\Gamma$, called the
\emph{inverse}, such that, for every $x\in\Gamma$,
$\bigl(x,\iota(x)\bigr)\in\Gamma_2$ and $\bigr(\iota(x),x\bigr)\in
\Gamma_2$, and
\begin{equation}
m\bigl(x,\iota(x)\bigr)=\varepsilon\bigl(\alpha(x)\bigr)\,,\quad
m\bigl(\iota(x),x\bigr)=\varepsilon\bigl(\beta(x)\bigr)\,.
\end{equation}
\end{description}
The sets $\Gamma$ and $\Gamma_0$ are called, respectively, the
\emph{total space} and the \emph{set of units} of the groupoid,
which is itself denoted by $\Gamma\toto^\alpha_\beta\Gamma_0$.
\end{defi}

\subsection{Identification and notations} In what follows, by
means of the injective map $\varepsilon$, we will identify the set
of units $\Gamma_0$ with se subset $\varepsilon(\Gamma_0)$ of
$\Gamma$. Therefore $\varepsilon$ will be the canonical injection
in $\Gamma$ of its subset $\Gamma_0$.
\par
For $x$ and $y\in\Gamma$, we will sometimes write $x.y$, or even
simply $xy$ for $m(x,y)$, and $x^{-1}$ for $\iota(x)$. Also we
will write \lq\lq the groupoid $\Gamma$\rq\rq\ for \lq\lq the
groupoid $\Gamma\toto^\alpha_\beta \Gamma_0$.

\subsection{Properties and comments}
The above definitions have the following consequences.

\subsubsection{Involutivity of the inverse map} The inverse
map $\iota$ is involutive:
 \begin{equation}
 \iota\circ\iota=\id_\Gamma\,.
 \end{equation}
 We have indeed, for any
$x\in\Gamma$,
 \[
 \begin{split}
 \iota\circ\iota(x)&=m\bigl(\iota\circ\iota(x), \beta\bigl(\iota\circ\iota(x)\bigr)\bigr)
        =m\bigl(\iota\circ\iota(x),\beta(x)\bigr)
        =m\bigl(\iota\circ\iota(x),m\bigl(\iota(x),x\bigr)\bigr)\\
       &=m\bigl(m\bigl(\iota\circ\iota(x),\iota(x)\bigr),x\bigr)
        =m\bigl(\alpha(x),x\bigr)
         =x\,.\\
 \end{split}
 \]

\subsubsection{Unicity of the inverse} Let $x$ and $y\in
\Gamma$ be such that
\[m(x,y)=\alpha(x)\quad\hbox{and}\quad m(y,x)=\beta(x)\,.\]
Then we have
 \[\begin{split}
   y&=m\bigl(y,\beta(y)\bigr)=m\bigl(y,\alpha(x)\bigr)
     =m\bigl(y,m\bigl(x,\iota(x)\bigr)\bigr)
     =m\bigl(m(y,x),\iota(x)\bigr)\\
    &=m\bigl(\beta(x),\iota(x)\bigr)
     =m\bigl(\alpha\bigl(\iota(x)\bigr),\iota(x)\bigr)
     =\iota(x)\,.\\
 \end{split}\]
Therefore for any $x\in \Gamma$, the unique $y\in\Gamma$ such that
$m(y,x)=\beta(x)$ and $m(x,y)=\alpha(x)$ is $\iota(x)$.

\subsubsection{The fibers of $\alpha$ and $\beta$ and the isotropy
groups}\label{isotropy} The target map $\alpha$ (resp. the source
map $\beta$) of a groupoid $\Gamma\toto^\alpha_\beta\Gamma_0$
determines an equivalence relation on $\Gamma$: two elements $x$
and $y\in\Gamma$ are said to be $\alpha$-equivalent (resp.
$\beta$-equivalent) if $\alpha(x)=\alpha(y)$ (resp. if
$\beta(x)=\beta(y)$). The corresponding equivalence classes are
called the \emph{$\alpha$-fibers} (resp. the
\emph{$\beta$-fibers}) of the groupoid. They are of the form
$\alpha^{-1}(u)$ (resp. $\beta^{-1}(u)$), with $u\in\Gamma_0$.
\par
For each unit $u\in\Gamma_0$, the subset
 \begin{equation}
 \Gamma_u=\alpha^{-1}(u)\cap\beta^{-1}(u)
 =\bigl\{\,x\in\Gamma;\alpha(x)=\beta(x)=u\,\bigr\}
 \end{equation}
is called the {\it isotropy group\/} of $u$. It is indeed a group,
with the restrictions of $m$ and $\iota$ as composition law and
inverse map.

\subsubsection{A way to visualize groupoids} We have seen
(Figure~1) a way in which groupoids may be visualized, by using
arrows for elements in $\Gamma$ and points for elements in
$\Gamma_0$. There is another, very useful way to visualize
groupoids, shown on Figure~2.

\begin{figure}[h]
\includegraphics{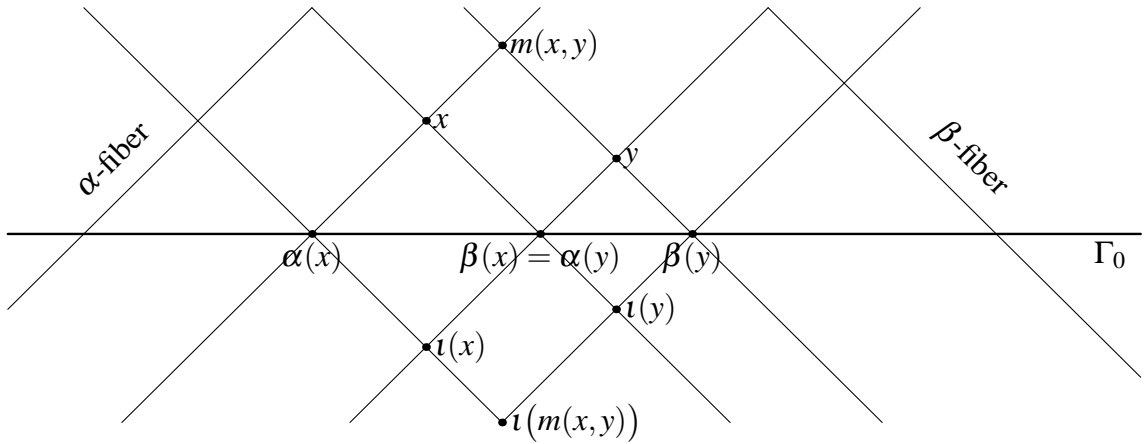}
\caption{A way to visualize groupoids.}
\end{figure}
The total space $\Gamma$ of the groupoid is represented as a
plane, and the set $\Gamma_0$ of units as a straight line in that
plane. The $\alpha$-fibers (resp. the $\beta$-fibers) are
represented as parallel straight lines, transverse to $\Gamma_0$.

\subsection{Examples of groupoids}\hfill

\subsubsection{The groupoid of pairs} Let $E$ be a set. The
{\it groupoid of pairs\/} of elements in $E$ has, as its total
space, the product space $E\times E$. The diagonal
$\Delta_E=\bigl\{\,(x,x);x\in E\,\bigr\}$ is its set of units, and
the target and source maps are
 $$\alpha:(x,y)\mapsto(x,x)\,,\quad \beta:(x,y)\mapsto(y,y)\,.$$
Its composition law $m$ and inverse map $\iota$ are
 $$m\bigl((x,y),(y,z)\bigr)=(x,z)\,,\quad
 \iota\bigl((x,y)\bigr)=(x,y)^{-1}=(y,x)\,.$$

\subsubsection{Groups} A group $G$ is a groupoid with set of units
$\{e\}$, with only one element $e$, the unit element of the group.
The target and source maps are both equal to the constant map
$x\mapsto e$.

\begin{defis}
A \emph{topological groupoid} is a groupoid
$\Gamma\toto^\alpha_\beta\Gamma_0$ for which $\Gamma$ is a (maybe
non Hausdorff) topological space, $\Gamma_0$ a Hausdorff
topological subspace of $\Gamma$, $\alpha$ and $\beta$ surjective
continuous maps, $m:\Gamma_2\to\Gamma$ a continuous map and
$\iota:\Gamma\to\Gamma$ an homeomorphism.
\par\smallskip
A \emph{Lie groupoid} is a groupoid
$\Gamma\toto^\alpha_\beta\Gamma_0$ for which $\Gamma$ is a smooth
(maybe non Hausdorff) manifold, $\Gamma_0$ a smooth Hausdorff
submanifold of $\Gamma$, $\alpha$ and $\beta$ smooth surjective
submersions (which implies that $\Gamma_2$ is a smooth submanifold
of $\Gamma\times\Gamma$), $m:\Gamma_2\to\Gamma$ a smooth map and
$\iota:\Gamma\to\Gamma$ a smooth diffeomorphism.
\end{defis}

\goodbreak
\subsection{Properties of Lie groupoids}\hfill
\nobreak
\subsubsection{Dimensions} Let $\Gamma\toto^\alpha_\beta\Gamma_0$ be a Lie
groupoid. Since $\alpha$ and $\beta$ are submersions, for any
$x\in\Gamma$, the $\alpha$-fiber
$\alpha^{-1}\bigl(\alpha(x)\bigr)$ and the $\beta$-fiber
$\beta^{-1}\bigl(\beta(x)\bigr)$ are submanifolds of $\Gamma$,
both of dimension $\dim\Gamma-\dim\Gamma_0$. The inverse map
$\iota$, restricted to the $\alpha$-fiber through $x$ (resp. the
$\beta$-fiber through $x$) is a diffeomorphism of that fiber onto
the $\beta$-fiber through $\iota(x)$ (resp. the $\alpha$-fiber
through $\iota(x)$). The dimension of the submanifold $\Gamma_2$
of composable pairs in $\Gamma\times \Gamma$ is
$2\dim\Gamma-\dim\Gamma_0$.

\subsubsection{The tangent bundle of a Lie groupoid}
Let $\Gamma\toto^\alpha_\beta\Gamma_0$ be a Lie
groupoid. Its tangent bundle $T\Gamma$ is a Lie groupoid, with
$T\Gamma_0$ as set of units, $T\alpha:T\Gamma\to T\Gamma_0$ and
$T\beta:T\Gamma\to T\Gamma_0$ as target and source maps. Let us
denote by $\Gamma_2$ the set of composable pairs in
$\Gamma\times\Gamma$, by $m:\Gamma_2\to\Gamma$ the composition law
and by $\iota:\Gamma\to\Gamma$ the inverse. Then the set of
composable pairs in $T\Gamma\times T\Gamma$ is simply $T\Gamma_2$,
the composition law on $T\Gamma$ is $Tm:T\Gamma_2\to T\Gamma$ and
the inverse is $T\iota:T\Gamma\to T\Gamma$.

When the groupoid $\Gamma$ is a Lie group $G$, the Lie groupoid
$TG$ is a Lie group too.

We will see below that the cotangent bundle of a Lie groupoid is a
Lie groupoid, and more precisely a symplectic groupoid.

\subsubsection{Isotropy groups} For each unit $u\in\Gamma_0$ of a Lie groupoid, the
isotropy group $\Gamma_u$ (defined in \ref{isotropy}) is a Lie
group.

\subsection{Examples of topological and Lie groupoids}\hfill

\subsubsection{Topological groups and Lie groups} A topological
group (resp. a Lie group) is a topological groupoid (resp. a Lie
groupoid) whose set of units has only one element $e$.

\subsubsection{Vector bundles} A smooth vector bundle $\pi:E\to M$
on a smooth manifold $M$ is a Lie groupoid, with the base $M$ as
set of units (identified with the image of the zero section); the
source and target maps both coincide with the projection $\pi$,
the product and the inverse maps are the addition $(x,y)\mapsto
x+y$ and the opposite map $x\mapsto -x$ in the fibers.

\subsubsection{The fundamental groupoid of a topological space} Let
$M$ be a topological space. A \emph{path} in $M$ is a continuous
map $\gamma:[0,1]\to M$. We denote by $[\gamma]$ the homotopy
class of a path $\gamma$ and by $\Pi(M)$ the set of homotopy
classes of paths in $M$ (with fixed endpoints). For $[\gamma]\in
\Pi(M)$, we set
 $\alpha\bigl([\gamma]\bigr)=\gamma(1)$,
   $\beta\bigl([\gamma]\bigr)=\gamma(0)$,
where $\gamma$ is any representative of the class $[\gamma]$. The
concatenation of paths determines a well defined composition law
on $\Pi(M)$, for which $\Pi(M)\toto^\alpha_\beta M$ is a
topological groupoid, called the \emph{fundamental groupoid} of
$M$. The inverse map is $[\gamma]\mapsto [\gamma^{-1}]$, where
$\gamma$ is any representative of $[\gamma]$ and $\gamma^{-1}$ is
the path $t\mapsto\gamma(1-t)$. The set of units is $M$, if we
identify a point in $M$ with the homotoppy class of the constant
path equal to that point.
\par\smallskip
When $M$ is a smooth manifold, the same construction can be made
with piecewise smooth paths, and the fundamental groupoid
$\Pi(M)\toto^\alpha_\beta M$ is a Lie groupoid.

\section{Symplectic and Poisson groupoids}

\subsection{Symplectic and Poisson geometry}

Let us recall some definitions and results in symplectic and
Poisson geometry, used in the next sections.

\subsubsection{Symplectic manifolds} A \emph{symplectic form}
on a smooth manifold $M$ is
a differential $2$-form $\omega$, which is closed, {\it i.e.}
which satisfies
\begin{equation}
d\omega=0\,,
\end{equation}
and nondegenerate, {\it i.e.} such that for each point $x\in M$ an
each nonzero vector $v\in T_xM$, there exists a vector $w\in T_xM$
such that $\omega(v,w)\neq 0$. Equipped with the symplectic form
$\omega$, a smooth manifold $M$ is called a \emph{symplectic
manifold} and denoted by $(M,\omega)$.

The dimension of a symplectic manifold is always even.

\subsubsection{The Liouville form on a cotangent bundle} Let $N$ be a smooth
manifold, and $T^*N$ be its
cotangent bundle. The \emph{Liouville form} on $T^*N$ is the
$1$-form $\theta$ such that, for any $\eta\in T^*N$ and $v\in
T_\eta(T^*N)$,
  \begin{equation}
  \theta(v)=\bigl\langle\eta,T\pi_N(v)\bigr\rangle\,,
  \end{equation}
where $\pi_N:T^*N\to N$ is the canonical projection.

The $2$-form $\omega=d\theta$ is symplectic, and is called the
\emph{canonical symplectic form} on the cotangent bundle $T^*N$.

\subsubsection{Poisson manifolds} A \emph{Poisson manifold} is a smooth manifold $P$
equipped with a bivector field (i.e. a smooth section of
$\bigwedge^2TP$) $\Pi$ which satisfies
 \begin{equation}
 [\Pi,\Pi]=0\,,
 \end{equation}
the bracket on the left hand side being the Schouten bracket. The
bivector field $\Pi$ will be called the \emph{Poisson structure}
on $P$. It allows us to define a composition law on the space
$C^\infty(P,\RR)$ of smooth functions on $P$, called the
\emph{Poisson bracket} and denoted by $(f,g)\mapsto\{f,g\}$, by
setting, for all $f$ and $g\in C^\infty(P,\RR)$ and $x\in P$,
\begin{equation}
 \{f,g\}(x)=\Pi\bigl(df(x),dg(x)\bigr)\,.
 \end{equation}
That composition law is skew-symmetric and satisfies the Jacobi
identity, therefore turns $C^\infty(P,\RR)$ into a Lie algebra.

\subsubsection{Hamiltonian vector fields} Let
$(P,\Pi)$ be a Poisson manifold. We denote by
$\Pi^\sharp: T^*P\to TP$ the vector bundle map defined by
\begin{equation}
 \bigl\langle\eta,\Pi^\sharp(\zeta)\bigr\rangle=\Pi(\zeta,\eta)\,,
 \end{equation}
where $\zeta$ and $\eta$ are two elements in the same fiber of
$T^*P$. Let $f:P\to \RR$ be a smooth function on $P$. The vector
field $X_f=\Pi^\sharp(df)$ is called the \emph{Hamiltonian vector
field} associated to $f$. If $g:P\to\RR$ is another smooth
function on $P$, the Poisson bracket $\{f,g\}$ can be written
 \begin{equation}
 \{f,g\}=\bigl\langle dg,\Pi^\sharp(df)\bigr\rangle
          =-\bigl\langle df,\Pi^\sharp(dg)\bigr\rangle\,.
 \end{equation}

\subsubsection{The canonical Poisson structure on a symplectic manifold}
Every symplectic manifold $(M,\omega)$ has a Poisson
structure, associated to its symplectic structure, for which the
vector bundle map $\Pi^\sharp:T^*M \to M$ is the inverse of the
vector bundle isomorphism $v\mapsto -i(v)\omega$. We will always
consider that a symplectic manifold is equipped with that Poisson
structure, unless otherwise specified.

\subsubsection{The KKS Poisson structure}
\label{KKS}
Let $\cal G$
be a finite-dimensional Lie algebra.
Its dual space ${\cal G}^*$ has a natural Poisson structure, for
which the bracket of two smooth functions $f$ and $g$ is
 \begin{equation}
 \{f,g\}(\xi)=\bigl\langle\xi,\bigl[df(\xi),dg(\xi)\bigr]\bigr\rangle\,,
 \end{equation}
with $\xi\in{\cal G}^*$, the differentials $df(\xi)$ and $dg(\xi)$
being considered as elements in $\cal G$, identified with its
bidual ${\cal G}^{**}$. It is called the KKS Poisson structure on
${\cal G}^*$ (for Kirillov, Kostant and Souriau).

\subsubsection{Poisson maps} Let $(P_1,\Pi_1)$ and $(P_2,\Pi_2)$ be two Poisson
manifolds. A smooth map $\varphi:P_1\to P_2$ is called a
\emph{Poisson map} if for every pair $(f,g)$ of smooth functions
on $P_2$,
 \begin{equation}
 \{\varphi^*f,\varphi^*g\}_1=\varphi^*\{f,g\}_2\,.
 \end{equation}

\subsubsection{Product Poisson structures} The product
$P_1\times P_2$ of two Poisson manifolds
$(P_1,\Pi_1)$ and $(P_2,\Pi_2)$ has a natural Poisson structure:
it is the unique Poisson structure for which the bracket of
functions of the form $(x_1,x_2)\mapsto f_1(x_1)f_2(x_2)$ and
$(x_1,x_2)\mapsto g_1(x_1)g_2(x_2)$, where $f_1$ and $g_1\in
C^\infty(P_1,\RR)$, $f_2$ and $g_2\in C^\infty(P_2,\RR)$, is
 $$(x_1,x_2)\mapsto\{f_1,g_1\}_1(x_1)\{f_2,g_2\}_2(x_2)\,.$$
The same property holds for the product of any finite number of
Poisson manifolds.

\subsubsection{Symplectic orthogonality} Let $(V,\omega)$
be a symplectic vector space, that
means a real, finite-dimensional vector space $V$ with a
skew-symmetric nondegenrate bilinear form $\omega$. Let $W$ be a
vector subspace of $V$. The \emph{symplectic orthogonal} of $W$ is
 \begin{equation}
 \orth W=\bigl\{\,v\in V;\omega(v,w)=0\ \hbox{for all}\ w\in
 W\,\}\,.
 \end{equation}
It is a vector subspace of $V$, which satisfies
 \[\dim W+\dim(\orth W)=\dim V\,,\quad \orth(\orth W)=W\,.\]
The vector subspace $W$ is said to be \emph{isotropic} if
$W\subset\orth W$, \emph{coisotropic} if $\orth W\subset W$ and
\emph{Lagrangian} if $W=\orth W$. In any symplectic vector space,
there are many Lagrangian subspaces, therefore the dimension of a
symplectic vector space is always even; if $\dim V=2n$, the
dimension of an isotropic (resp. coisotropic, resp. Lagrangian)
vector subspace is $\leq n$ (resp. $\geq n$, resp. $=n$).

\subsubsection{Coisotropic and Lagrangian submanifolds} A submanifold
$N$ of a Poisson manifold $(P,\Pi)$ is
said to be \emph{coisotropic} if the bracket of two smooth
functions, defined on an open subset of $P$ and which vanish on
$N$, vanishes on $N$ too. A submanifold $N$ of a symplectic
manifold $(M,\omega)$ is coisotropic if and only if for each point
$x\in N$, the vector subspace $T_xN$ of the symplectic vector
space $\bigl(T_xM,\omega(x)\bigr)$ is coisotropic. Therefore, the
dimension of a coisotropic submanifold in a $2n$-dimensional
symplectic manifold is $\geq n$; when it is equal to $n$, the
submanifold $N$ is said to be \emph{Lagrangian}.

\subsubsection{Poisson quotients} Let $\varphi:M\to P$ be a
surjective submersion of a
symplectic manifold $(M,\omega)$ onto a manifold $P$. The manifold
$P$ has a Poisson structure $\Pi$ for which $\varphi$ is a Poisson
map if and only if $\orth(\ker T\varphi)$ is integrable. When that
condition is satisfied, that Poisson structure on $P$ is unique.

\subsubsection{Poisson Lie groups} A \emph{Poisson Lie group}
is a Lie group $G$ with a Poisson structure $\Pi$, such that the
product $(x,y)\mapsto xy$ is a Poisson map from $G\times G$,
endowed with the product Poisson structure, into $(G,\Pi)$. The
Poisson structure of a Poisson Lie group $(G,\Pi)$ always vanishes
at the unit element $e$ of $G$. Therefore the Poisson structure of
a Poisson Lie group never comes from a symplectic structure on
that group.

\begin{defis}
A \emph{symplectic groupoid} (resp. a \emph{Poisson groupoid}) is
a Lie groupoid $\Gamma\toto^\alpha_\beta\Gamma_0$ with a
symplectic form $\omega$ on $\Gamma$ (resp. with a Poisson
structure $\Pi$ on $\Gamma$) such that the graph of the
composition law $m$
 \[\bigl\{\,(x,y,z)\in\Gamma\times\Gamma\times\Gamma;
 (x,y)\in\Gamma_2\ \hbox{and}\ z=m(x,y)\,\bigr\}
 \]
is a Lagrangian submanifold (resp. a coisotropic submanifold) of
$\Gamma\times\Gamma\times\overline\Gamma$ with the product
symplectic form (resp. the product Poisson structure), the first
two factors $\Gamma$ being endowed with the symplectic form
$\omega$ (resp. with the Poisson structure $\Pi$), and the third
factor $\overline\Gamma$ being $\Gamma$ with the symplectic form
$-\omega$ (resp. with the Poisson structure $-\Pi$).
\end{defis}

The next theorem states important properties of symplectic and
Poisson groupoids.

\begin{theo}
\label{symplgroup} Let $\Gamma\toto^\alpha_\beta\Gamma_0$ be a
symplectic groupoid with symplectic $2$-form $\omega$ (resp. a
Poisson groupoid with Poisson structure $\Pi$). We have the
following properties.

\smallskip\noindent{\rm 1.} For a symplectic groupoid, given
any point $c \in\Gamma$, each one of the two vector subspaces of
the symplectic vector space $\bigl(T_c\Gamma, \omega(c)\bigr)$,
$T_c\bigl(\beta^{-1}\bigl(\beta(c)\bigr)\bigr)$ and
$T_c\bigl(\alpha^{-1}\bigl(\alpha(c)\bigr)\bigr)$, is the
symplectic orthogonal of the other one. For a symplectic or
Poisson groupoid, if $f$ is a smooth function whose restriction to
each $\alpha$-fiber is constant, and $g$ a smooth function whose
restriction to each $\beta$-fiber is constant, then the Poisson
bracket $\{f,g\}$ vanishes identically.

\smallskip\noindent{\rm 2.} The submanifold of units $\Gamma_0$
is a Lagrangian submanifold of the symplectic manifold
$(\Gamma,\omega)$ (resp. a coisotropic submanifold of the Poisson
manifold $(\Gamma,\Pi)$).

\smallskip\noindent{\rm 3.} The inverse map $\iota:\Gamma\to\Gamma$ is an
antisymplectomorphism of $(\Gamma,\omega)$, i.e. it satisfies
$\iota^*\omega=-\omega$ (resp an anti-Poisson diffeomerphism of
$(\Gamma,\Pi)$, i.e. it satisfies $\iota_*\Pi=-\Pi$).
\end{theo}

\begin{coro}
\label{poissonbase}
Let $\Gamma\toto^\alpha_\beta\Gamma_0$ be a
symplectic groupoid with symplectic $2$-form $\omega$ (resp. a
Poisson groupoid with Poisson structure $\Pi$). There exists on
$\Gamma_0$ a unique Poisson structure $\Pi_0$ for which
$\alpha:\Gamma\to\Gamma_0$ is a Poisson map, and
$\beta:\Gamma\to\Gamma_0$ an anti-Poisson map ({\rm i.e.} $\beta$
is a Poisson map when $\Gamma_0$ is equipped with the Poisson
structure $-\Pi_0$).
\end{coro}

\subsection{Examples of symplectic and Poisson groupoids}
\subsubsection{The cotangent bundle of a Lie groupoid}
Let $\Gamma\toto^\alpha_\beta\Gamma_0$ be a Lie groupoid.
\goodbreak
We have seen above that its tangent bundle $T\Gamma$ has a Lie
groupoid structure, determined by that of $\Gamma$. Similarly (but
much less obviously) the cotangent bundle $T^*\Gamma$ has a Lie
groupoid structure determined by that of $\Gamma$. The set of
units is the conormal bundle to the submanifold $\Gamma_0$ of
$\Gamma$, denoted by ${\cal N}^*\Gamma_0$. We recall that ${\cal
N}^*\Gamma_0$ is the vector sub-bundle of $T^*_{\Gamma_0}\Gamma$
(the restriction to $\Gamma_0$ of the cotangent bundle
$T^*\Gamma$) whose fiber ${\cal N}^*_p\Gamma_0$ at a point
$p\in\Gamma_0$ is
 \[{\cal N}^*_p\Gamma_0=\bigl\{\,\eta\in
 T^*_p\Gamma;\langle\eta,v\rangle=0\ \hbox{for all}\ v\in
 T_p\Gamma_0\,\bigr\}\,.\]
To define the target and source maps of the Lie algebroid
$T^*\Gamma$, we introduce the notion  of \emph{bisection} through
a point $x\in \Gamma$. A bisection through $x$ is a submanifold
$A$ of $\Gamma$, with $x\in A$, transverse both to the
$\alpha$-fibers and to the $\beta$-fibers, such that the maps
$\alpha$ and $\beta$, when restricted to $A$, are diffeomorphisms
of $A$ onto open subsets $\alpha(A)$ and $\beta(A)$ of $\Gamma_0$,
respectively. For any point $x\in M$, there exist bisections
through $x$. A bisection $A$ allows us to define two smooth
diffeomorphisms between open subsets of $\Gamma$, denoted by $L_A$
and $R_A$ and called the \emph{left} and \emph{right translations
by} $A$, respectively. They are defined by
 \[L_A:\alpha^{-1}\bigl(\beta(A)\bigr)\to\alpha^{-1}
 \bigl(\alpha(A)\bigr)\,,\quad
 L_A(y)=m\bigl(\beta|_A^{-1}\circ\alpha(y),y\bigr)\,,\]
and
\[R_A:\beta^{-1}\bigl(\alpha(A)\bigr)\to\beta^{-1}
 \bigl(\beta(A)\bigr)\,,\quad
 R_A(y)=m\bigl(y, \alpha|_A^{-1}\circ\beta(y)\bigr)\,.\]
The definitions of the target and source maps for $T^*\Gamma$ rest
on the following properties. Let $x$ be a point in $\Gamma$ and
$A$ be a bisection through $x$. The two vector subspaces,
$T_{\alpha(x)}\Gamma_0$ and $\ker T_{\alpha(x)}\beta$, are
complementary in $T_{\alpha(x)}\Gamma$. For any $v\in
T_{\alpha(x)}\Gamma$, $v-T\beta(v)$ is in $\ker
T_{\alpha(x)}\beta$. Moreover, $R_A$ maps the fiber
$\beta^{-1}\bigl(\alpha(x)\bigr)$ onto the fiber
$\beta^{-1}\bigl(\beta(x)\bigr)$, and its restriction to that
fiber does not depend on the choice of $A$; its depends only on
$x$. Therefore $TR_A\bigl(v-T\beta(v)\bigr)$ is in $\ker T_x\beta$
and does not depend on  the choice of $A$. We can define the map
$\widehat\alpha$ by setting, for any $\xi\in T^*_x\Gamma$ and any
$v\in T_{\alpha(x)}\Gamma$,
 \[\bigl\langle\widehat\alpha(\xi),v\bigr\rangle=\bigl\langle
 \xi, TR_A\bigl(v-T\beta(v)\bigr)\bigr\rangle\,.\]
Similarly, we define $\widehat\beta$ by setting, for any $\xi\in
T^*_x\Gamma$ and any $w\in T_{\beta(x)}\Gamma$,
 \[\bigl\langle\widehat\beta(\xi),w\bigr\rangle=\bigl\langle
 \xi, TL_A\bigl(w-T\alpha(w)\bigr)\bigr\rangle\,.\]
We see that $\widehat\alpha$ and $\widehat\beta$ are unambiguously
defined, smooth and take their values in the submanifold ${\cal
N}^*\Gamma_0$ of $T^*\Gamma$. They satisfy
 \[\pi_\Gamma\circ\widehat\alpha=\alpha\circ\pi_\Gamma\,,\quad
   \pi_\Gamma\circ\widehat\beta=\beta\circ\pi_\Gamma\,,\]
where $\pi_\Gamma:T^*\Gamma\to\Gamma$ is the cotangent bundle
projection.

Let us now define the composition law $\widehat m$ on $T^*\Gamma$.
Let $\xi\in T^*_x\Gamma$ and $\eta\in T^*_y\Gamma$ be such that
$\widehat\beta(\xi)=\widehat\alpha(\eta)$. That implies
$\beta(x)=\alpha(y)$. Let $A$ be a bisection through $x$ and $B$ a
bisection through $y$. There exist a unique $\xi_{h\alpha}\in
T^*_{\alpha(x)}\Gamma_0$ and a unique $\eta_{h\beta}\in
T^*_{\beta(y)}\Gamma_0$ such that
 \[\xi=(L_A^{-1})^*\bigl(\widehat\beta(\xi)\bigr)+\alpha^*_x\xi_{h\alpha}\,,\quad
   \eta=(R_B^{-1})^*\bigl(\widehat\alpha(\xi)\bigr)+\beta^*_y\eta_{h\beta}\,.\]
Then $\widehat m(\xi,\eta)$ is given by
 \[\widehat
 m(\xi,\eta)=\alpha^*_{xy}\xi_{h\alpha}+\beta^*_{xy}\eta_{h\beta}
 +(R_B^{-1})^*(L_A^{-1})^*\bigl(\widehat\beta(x)\bigr)\,.\]
We observe that in the last term of the above expression we can
replace $\widehat\beta(\xi)$ by $\widehat\alpha(\eta)$, since
these two expressions are equal, and that
$(R_B^{-1})^*(L_A^{-1})^*=(L_A^{-1})^*(R_B^{-1})^*$, since $R_B$
and $L_A$ commute.

Finally, the inverse $\widehat \iota$ in $T^*\Gamma$ is $\iota^*$.

With its canonical symplectic form,
$T^*\Gamma\toto^{\widehat\alpha}_{\widehat\beta}{\cal
N}^*\Gamma_0$ is a symplectic groupoid.

 When the Lie groupoid
$\Gamma$ is a Lie group $G$, the Lie groupoid $T^*G$ is not a Lie
group, contrary to what happens for $TG$. This shows that the
introduction of Lie groupoids is not at all artificial: when
dealing with Lie groups, Lie groupoids are already with us! The
set of units of the Lie groupoid $T^*G$ can be identified with
${\cal G}^*$ (the dual of the Lie algebra $\cal G$ of $G$),
identified itself with $T^*_eG$ (the cotangent space to $G$ at the
unit element $e$). The target map $\widehat\alpha:T^*G\to T^*_eG$
(resp. the source map $\widehat\beta:T^*G\to T^*_eG$) associates
to each $g\in G$ and $\xi\in T^*_gG$, the value at the unit
element $e$ of the right-invariant $1$-form (resp., the
left-invariant $1$-form) whose value at $x$ is $\xi$.

\subsubsection{Poisson Lie groups as Poisson groupoids}
Poisson groupoids were introduced by Alan~Weinstein as a
generalization of both symplectic groupoids and Poisson Lie
groups. Indeed, a Poisson Lie group is a Poisson groupoid with a
set of units reduced to a single element.

\section{Lie algebroids}

The notion of a Lie algebroid, due to Jean~Pradines, is related to
that of a Lie groupoid in the same way as the notion of a Lie
algebra is related to that of a Lie group.

\begin{defi} A \emph{Lie algebroid} over a smooth manifold $M$ is a
smooth vector bundle $\pi:A\to M$ with base $M$, equipped with
\begin{description}
\item{--}  a composition law $(s_1,s_2)\mapsto\{s_1,s_2\}$ on the
space
  $\Gamma^\infty(\pi)$ of smooth sections of $\pi$, called the \emph{bracket},
  for which that space is a Lie algebra,

\item{--}  a vector bundle map $\rho:A\to TM$, over the identity
map of $M$, called the \emph{anchor map}, such that, for all $s_1$
and  $s_2\in\Gamma^\infty(\pi)$ and all $f\in C^\infty(M,\RR)$,
  \begin{equation}
  \{s_1,fs_2\}=f\{s_1,s_2\}+\bigl((\rho\circ
  s_1).f\bigr)s_2\,.
  \end{equation}
\end{description}
\end{defi}

\subsection{Examples}\hfill

\subsubsection{Lie algebras}
A finite-dimensional Lie algebra is a Lie algebroid (with a base
reduced to a point and the zero map as anchor map).

\subsubsection{Tangent bundles and their integrable sub-bundless}
\label{tangbundasliealgebroid}
A tangent bundle  $\tau_M:TM\to M$
to a smooth manifold $M$ is a Lie algebroid, with the usual
bracket of vector fields on $M$ as composition law, and the
identity map as anchor map. More generally, any integrable vector
sub-bundle $F$ of a tangent bundle $\tau_M:TM\to M$ is a Lie
algebroid, still with the bracket of vector fields on $M$ with
values in $F$ as composition law and the canonical injection of
$F$ into $TM$ as anchor map.

\subsubsection{The cotangent bundle of a Poisson manifold}
\label{bracketofforms}
Let $(P,\Pi)$ be a Poisson manifold. Its
cotangent bundle $\pi_P:T^*P\to P$ has a Lie algebroid structure,
with $\Pi^\sharp:T^*P\to TP$ as anchor map. The composition law is
the \emph{bracket of $1$-forms}. It will be denoted by
$(\eta,\zeta)\mapsto[\eta,\zeta]$ (in order to avoid any confusion
with the Poisson bracket of functions). It is given by the
formula, in which $\eta$ and $\zeta$ are $1$-forms and $X$ a
vector field on $P$,
 \begin{equation}
 \bigl\langle[\eta,\zeta],X\bigr\rangle=
 \Pi\bigl(\eta,d\langle\zeta,X\rangle\bigr)
 +\Pi\bigl(d\langle\eta,X\rangle,\zeta\bigr)
 +\bigl({\cal L}(X)\Pi\bigr)(\eta,\zeta)\,.
 \end{equation}
We have denoted by ${\cal L}(X)\Pi$ the Lie derivative of the
Poisson structure $\Pi$ with respect to the vector field $X$.
Another equivalent formula for that composition law is
 \begin{equation}
 [\zeta,\eta]={\cal L}(\Pi^\sharp\zeta)\eta
                 -{\cal L}(\Pi^\sharp\eta)\zeta
                 -d\bigl(\Pi(\zeta,\eta)\bigr)\,.
 \end{equation}
The bracket of $1$-forms is related to the Poisson bracket of
functions by
 \begin{equation}
 [df,dg]=d\{f,g\}\quad\hbox{for all}\ f\ \hbox{and}\ g\in
 C^\infty(P,\RR)\,.
 \end{equation}

\subsection{Properties of Lie algebroids} Let $\pi:A\to M$ be a
Lie algebroid with anchor map $\rho:A\to TM$.

\subsubsection{A Lie algebras homomorphism}
For any pair $(s_1,s_2)$ of smooth sections of
$\pi$,
 \[\rho\circ\{s_1,s_2\}=[\rho\circ s_1,\rho\circ s_2]\,,\]
which means that the map $s\mapsto\rho\circ s$ is a Lie algebra
homomorphism from the Lie algebra of smooth sections of $\pi$ into
the Lie algebra of smooth vector fields on $M$.

\subsubsection{The generalized Schouten bracket}
\label{genschoubrack}
The composition law
$(s_1,s_2)\mapsto\{s_1,s_2\}$ on the space of sections of $\pi$
extends into a composition law on the space of sections of
exterior powers of $(A,\pi,M)$, which will be called the
\emph{generalized Schouten bracket}. Its properties are the same
as those of the usual Schouten bracket. When the Lie algebroid is
a tangent bundle $\tau_M:TM\to M$, that composition law reduces to
the usual Schouten bracket . When the Lie algebroid is the
cotangent bundle $\pi_P:T^*P\to P$ to a Poisson manifold
$(P,\Pi)$, the generalized Schouten bracket is the \emph{bracket
of forms} of all degrees on the Poisson manifold $P$, introduced
by J.-L.~Koszul, which extends the bracket of $1$-forms used in
\ref{bracketofforms}.

\subsubsection{The dual bundle of a Lie algebroid}
\label{genextdiff}
Let $\varpi:A^*\to M$ be the dual bundle of the
Lie algebroid $\pi:A\to M$. There exists on the space of sections
of its exterior powers a graded endomorphism $d_\rho$, of degree
$1$ (that means that if $\eta$ is a section of $\bigwedge^kA^*$,
$d_\rho(\eta)$ is a section of $\bigwedge^{k+1}A^*$). That
endomorphism satisfies
 \[d_\rho\circ d_\rho=0\,,\]
and its properties are essentially the same as those of the
exterior derivative of differential forms. When the Lie algebroid
is a tangent bundle $\tau_M:TM\to M$, $d_\rho$ is the usual
exterior derivative of differential forms.

We can develop on the spaces of sections of the exterior powers of
a Lie algebroid and of its dual bundle a differential calculus
very similar to the usual differential calculus of vector and
multivector fields and differential forms on a manifold. Operators
such as the interior product, the exterior derivative and the Lie
derivative can still be defined and have properties similar to
those of the corresponding operators for vector and multivector
fields and differential forms on a manifold.

The total space $A^*$ of the dual bundle of a Lie algebroid
$\pi:A\to M$ has a natural Poisson structure: a smooth section $s$
of $\pi$ can be considered as a smooth real-valued function on
$A^*$ whose restrictiion to each fiber $\varpi^{-1}(x)$ ($x\in M$)
is linear; that property allows us to extend the bracket of
sections of $\pi$ (defined by the Lie algebroid structure) to
obtain a Poisson bracket of functions on $A^*$. When the Lie
algebroid $A$ is a finite-dimensional Lie algebra $\cal G$, the
Poisson structure on its dual space ${\cal G}^*$ is the KKS
Poisson structure discussed in \ref{KKS}.

\subsection{The Lie algebroid of a Lie groupoid}

Let $\Gamma\toto^\alpha_\beta\Gamma_0$ be a Lie groupoid. Let
$A(\Gamma)$ be the intersection of $\ker T\alpha$ and
$T_{\Gamma_0}\Gamma$ (the tangent bundle $T\Gamma$ restricted to
the submanifold $\Gamma_0$). We see that $A(\Gamma)$ is the total
space of a vector bundle $\pi:A(\Gamma)\to\Gamma_0$, with base
$\Gamma_0$, the canonical projection $\pi$ being the map which
associates a point $u\in \Gamma_0$ to every vector in $\ker
T_u\alpha$. We will define a composition law on the set of smooth
sections of that bundle, and a vector bundle map
$\rho:A(\Gamma)\to T\Gamma_0$, for which
$\pi:A(\Gamma)\to\Gamma_0$ is a Lie algebroid, called the
\emph{Lie algebroid} of  the Lie groupoid
$\Gamma\toto^\alpha_\beta\Gamma_0$.

We observe first that for any point $u\in\Gamma_0$ and any point
$x\in\beta^{-1}(u)$, the map $L_x:y\mapsto L_xy=m(x,y)$ is defined
on the $\alpha$-fiber $\alpha^{-1}(u)$, and maps that fiber onto
the $\alpha$-fiber $\alpha^{-1}\bigl(\alpha(x)\bigr)$. Therefore
$T_uL_x$ maps the vector space $A_u=\ker T_u\alpha$ onto the
vector space $\ker T_x\alpha$, tangent at $x$ to the
$\alpha$-fiber $\alpha^{-1}\bigl(\alpha(x)\bigr)$. Any vector
$w\in A_u$ can therefore be extended into the vector field along
$\beta^{-1}(u)$, $x\mapsto \widehat w(x)=T_uL_x(w)$. More
generally, let $w:U\to A(\Gamma)$ be a smooth section of the
vector bundle $\pi:A(\Gamma)\to\Gamma_0$, defined on an open
subset $U$ of $\Gamma_0$. By using the above described
construction for every point $u\in U$, we can extend the section
$w$ into a smooth vector field  $\widehat w$, defined on the open
subset $\beta^{-1}(U)$ of $\Gamma$, by setting, for all $u\in U$
and $x\in \beta^{-1}(u)$,
 \[\widehat w(x)=T_uL_x\bigl(w(u)\bigr)\,.\]
We have defined an injective map $w\mapsto\widehat w$ from the
space of smooth local sections of $\pi:A(\Gamma)\to\Gamma_0$, onto
a subspace of the space of smooth vector fields defined on open
subsets of $\Gamma$. The image of that map is the space of smooth
vector fields $\widehat w$, defined on open subsets $\widehat U$
of $\Gamma$ of the form $\widehat U=\beta^{-1}(U)$, where $U$ is
an open subset of $\Gamma_0$, which satisfy the two properties:

\smallskip
\begin{description}
\item{(i)} $T\alpha\circ\widehat w=0$,

\smallskip
\item{(ii)} for every $x$ and $y\in\widehat U$ such that
$\beta(x)=\alpha(y)$, $T_yL_x\bigl(\widehat w(y)\bigr)=\widehat
w(xy)$.
\end{description}
\smallskip
These vector fields are called \emph{left invariant vector fields}
on $\Gamma$.

The space of left invariant vector fields on $\Gamma$ is closed
under the bracket operation.  We can therefore define a
composition law $(w_1,w_2)\mapsto\{w_1,w_2\}$ on the space of
smooth sections of the bundle $\pi:A(\Gamma)\to\Gamma_0$ by
defining $\{w_1,w_2\}$ as the unique section such that
 \[\widehat{\{w_1,w_2\}}=[\widehat w_1,\widehat w_2]\,.\]
Finally, we define the anchor map $\rho$ as the map $T\beta$
restricted to $A(\Gamma)$. With that composition law and that
anchor map, the vector bundle $\pi:A(\Gamma)\to\Gamma_0$ is a Lie
algebroid, called the \emph{Lie algebroid of} the Lie groupoid
$\Gamma\toto^\alpha_\beta\Gamma_0$.

We could exchange the roles of $\alpha$ and $\beta$ and use
\emph{right invariant} vector fields instead of left invariant
vector fields. The Lie algebroid obtained remains the same, up to
an isomorphism.

When the Lie groupoid $\Gamma\toto^\alpha_\beta$ is a Lie group,
its Lie algebroid is simply its Lie algebra.

\subsection{The Lie algebroid of a symplectic groupoid}
\label{liealgsympgroupoid}
 Let $\Gamma\toto^\alpha_\beta\Gamma_0$ be a symplectic
groupoid, with symplectic form $\omega$. As we have seen above,
its Lie algebroid $\pi:A\to \Gamma_0$ is the vector bundle whose
fiber, over each point $u\in \Gamma_0$, is $\ker T_u\alpha$. We
define a linear map $\omega_u^\flat:\ker T_u\alpha\to
T^*_u\Gamma_0$ by setting, for each $w\in \ker T_u\alpha$ and
$v\in T_u\Gamma_0$,
 \[\bigl\langle \omega_u^\flat(w),
 v\bigr\rangle=\omega_u(v,w)\,.\]
Since $T_u\Gamma_0$ is Lagrangian and $\ker T_u\alpha$
complementary to $T_u\Gamma_0$ in the symplectic vector space
$\bigl(T_u\Gamma,\omega(u)\bigr)$, the map $\omega_u^\flat$ is an
isomorphism from $\ker T_u\alpha$ onto $T^*_u\Gamma_0$. By using
that isomorphism for each $u\in\Gamma_0$, we obtain  a vector
bundle isomorphism of the Lie algebroid $\pi:A\to\Gamma_0$ onto
the cotangent bundle $\pi_{\Gamma_0}:T^*\Gamma_0\to\Gamma_0$.

As seen in Corollary~\ref{poissonbase}, the submanifold of units
$\Gamma_0$ has a unique Poisson structure $\Pi$ for which
$\alpha:\Gamma\to\Gamma_0$ is a Poisson map. Therefore, as seen in
\ref{bracketofforms}, the cotangent bundle
$\pi_{\Gamma_0}:T^*\Gamma_0\to\Gamma_0$ to the Poisson manifold
$(\Gamma_0,\Pi)$ has a Lie algebroid structure, with the bracket
of $1$-forms as composition law. That structure is the same as the
structure obtained as a direct image of the Lie algebroid
structure of $\pi:A(\Gamma)\to\Gamma_0$, by the above defined
vector bundle isomorphism of $\pi:A\to\Gamma_0$ onto the cotangent
bundle $\pi_{\Gamma_0}:T^*\Gamma_0\to\Gamma_0$. The Lie algebroid
of the symplectic groupoid $\Gamma\toto^\alpha_\beta\Gamma_0$ can
therefore be identified with the Lie algebroid
$\pi_{\Gamma_0}:T^*\Gamma_0\to\Gamma_0$, with its Lie algebroid
structure of cotangent bundle to the Poisson manifold
$(\Gamma_0,\Pi)$.

\subsection{The Lie algebroid of a Poisson groupoid}

The Lie algebroid $\pi:A(\Gamma)\to \Gamma_0$ of a Poisson
groupoid has an additional structure: its dual bundle
$\varpi:A(\Gamma)^*\to\Gamma_0$ also has a Lie algebroid
structure, compatible in a certain sense (indicated below) with
that of $\pi:A(\Gamma)\to\Gamma_0$ (K.~Mackenzie and P.~Xu,
Y.~Kosmann-Schwarzbach, Z.-J.~Liu and P.~Xu).

The compatibility condition between the two Lie algebroid
structures on the two vector bundles in duality $\pi:A\to M$ and
$\varpi:A^*\to M$ can be written as follows:
 \begin{equation}
 d_*[X,Y]={\cal L}(X)d_*Y-{\cal L}(Y)d_*X\,,
 \end{equation}
where $X$ and $Y$ are two sections of $\pi$, or, using the
generalized Schouten bracket (\ref{genschoubrack}) of sections of
exterior powers of the Lie algebroid  $\pi:A\to M$,
 \begin{equation}
 d_*[X,Y]=[d_*X,Y]+[X,d_*Y]\,.
 \end{equation}
In these formulae $d_*$ is the generalized exterior derivative,
which acts on the space of sections of exterior powers of the
bundle $\pi:A\to M$, considered as the dual bundle of the Lie
algebroid $\varpi:A^*\to M$, defined in \ref{genextdiff}.

These conditions are equivalent to the similar conditions obtained
by exchange of the roles of $A$ and $A^*$.

When the Poisson groupoid $\Gamma\toto^\alpha_\beta\Gamma_0$ is a
symplectic groupoid, we have seen (\ref{liealgsympgroupoid}) that
its Lie algebroid is the cotangent bundle
$\pi_{\Gamma_0}:T^*\Gamma_0\to\Gamma_0$ to the Poisson manifold
$\Gamma_0$ (equipped with the Poisson structure for which $\alpha$
is a Poisson map). The dual bundle is the tangent bundle
$\tau_{\Gamma_0}:T\Gamma_0\to\Gamma_0$, with its natural Lie
algebroid structure defined in \ref{tangbundasliealgebroid}.

When the Poisson groupoid is a Poisson Lie group $(G,\Pi)$, its
Lie algebroid is its Lie algebra $\cal G$. Its dual space ${\cal
G}$ has a Lie algebra structure, compatible with that of $\cal G$
in the above defined sense, and the pair $({\cal G},{\cal G}^*)$
is called a \emph{Lie bialgebra}.

Conversely, if the Lie algebroid of a Lie groupoid is a Lie
bialgebroid (that means, if there exists on the dual vector bundle
of that Lie algebroid a compatible structure of Lie algebroid, in
the above defined sense), that Lie groupoid has a Poisson
structure for which it is a Poisson groupoid (K.~Mackenzie and
P.~Xu).

\subsection{Integration of Lie algebroids}

According to Lie's third theorem, for any given finite-dimensional
Lie algebra, there exists a Lie group whose Lie algebra is
isomorphic to that Lie algebra. The same property is not true for
Lie algebroids and Lie groupoids. The problem of finding necessary
and sufficient conditions under which a given Lie algebroid is
isomorphic to the Lie algebroid of a Lie groupoid remained open
for more than 30 years, although partial results were obtained. A
complete solution of that problem was recently obtained by
M.~Crainic and R.L.~Fernandes. Let us briefly sketch their
results.
\par\smallskip
Let $\pi:A\to M$ be a Lie algebroid and $\rho:A\to TM$ its anchor
map. A smooth path $a:I=[0,1]\to A$ is said to be
\emph{admissible} if, for all $t\in I$, $\rho\circ
a(t)=\frac{d}{dt}(\pi\circ a)(t)$. When  the Lie algebroid $A$ is
the Lie algebroid of a Lie groupoid $\Gamma$, it can be shown that
each admissible path in $A$ is, in a natural way, associated to a
smooth path in $\Gamma$ starting from a unit and contained in an
$\alpha$-fiber. When we do not know whether $A$ is the Lie
algebroid of a Lie goupoid or not, the space of admissible paths
in $A$ still can be used to define a topological groupoid ${\cal
G}(A)$ with connected and simply connected $\alpha$-fibers, called
the \emph{Weinstein groupoid\/} of $A$. When ${\cal G}(A)$ is a
Lie groupoid, its Lie algebroid is isomorphic to $A$, and when $A$
is the Lie algebroid of a Lie groupoid $\Gamma$, ${\cal G}(A)$ is
a Lie groupoid and is the unique (up to an isomorphism) Lie
groupoid with connected and simply connected $\alpha$-fibers with
$A$ as Lie algebroid; moreover, ${\cal G}(A)$ is a covering
groupoid of an open sub-groupoid of $\Gamma$. Crainic and
Fernandes have obtained computable necessary and sufficient
conditions under which the topological groupoid ${\cal G}(A)$ is a
Lie groupoid, {\it i.e.} necessary and sufficient conditions under
which $A$ is the Lie algebroid of a Lie groupoid.

\section*{Key words}

Groupoids, Lie groupoids, Lie algebroids, symplectic groupoids,
Poisson groupoids, Poisson Lie groups, bisections.

\section*{Further reading}

The reader will find more about groupoids in the very nice survey
paper \cite{Wein5}, and in the books \cite{CannW}, \cite{Mack}.
More information about symplectic and Poisson geometry can be
found in \cite{LibMa}, \cite{Vai}, \cite{ortega-ratiu}. The
Schouten bracket is discussed in \cite{Vai}. The paper \cite{Xu3}
presents many properties of Poisson groupoids. Integration of Lie
algebroids is fully discussed in \cite{crainicfernandes1}. Papers
\cite{Kara} and \cite{Zak} introduced Lie groupoids for
applications to quantization. The books \cite{sudrhoda} and
\cite{MarsdenRatiu} are symposia proceedings which contain several
papers about Lie, symplectic and Poisson groupoids, and many
references.

\end{document}